\documentclass[11pt]{article}
\usepackage[colorlinks=true,breaklinks=true,citecolor=red,linkcolor=blue,urlcolor=red,unicode]{hyperref}
\usepackage{amsmath,amssymb}

\newcommand{\keywords}[1]{\footnote{\hspace{-7mm}
{\bf\scriptsize Keywords:} {\scriptsize \it #1}}}
\newcommand{\address}[1]{\footnote{\hspace{-7mm}
{\bf\scriptsize Corresponding author address:} {\bf\scriptsize #1}}}
\makeatletter
\newtheorem{theorem}{Theorem}[section]
\newtheorem{lemma}{Lemma}[section]
\newtheorem{corollary}{Corollary}[section]

\newtheorem{definition}{Definition}[section]
\newtheorem{remark}{\textit{Remark}}[section]
\newtheorem{example}{\textit{Example}}[section]
\numberwithin{equation}{section}

\def\proof{\noindent{\it Proof.}\hspace{3mm}}

\begin{document}
\title{\textbf{Normal Structure and Contractions which diminish Radius in Metric and Banach Spaces}}
\author{\textsc{Abdelkader Dehici, Najeh Redjel and Sami Atailia}}
\date{}
\maketitle
\begin{abstract} In this paper, by using admissible sets, we give some fixed point results for orbitally contractions which diminish the radius of invariant convex subsets and orbits. Furthermore, a characterization of the weak normal structure by the fixed point property associated with this class of mappings is established.
\end{abstract}


\address{Laboratory of Informatics and Mathematics\\
University of Souk-Ahras,P.O.Box 1553, Souk-Ahras 41000, Algeria\\
\texttt{E-mails:  dehicikader@yahoo.fr, najehredje@yahoo.fr, professoratailia@gmail.com}}

\section {Introduction}

In 1922, S. Banach has established his famous contraction principle in the case of complete metric spaces. The advantage of this result is that it assures the existence and uniqueness of the fixed point for contraction mappings, in addition to the convergence of their associated Picard sequences to the fixed point. However, the situation is complicated when we deal with nonexpensive mappings (the case where the contraction constant is equal to 1). To see this, it suffices to take the discrete metric space $(X, d_0)$ where $X = \{0, 1\}$ and $T: X \rightarrow X$ such that $T0 = 1$ and $T1 = 0$. Obviously, $T$ is nonexpansive since $d_0(Tx,Ty) = d_0(x,y)$ for all $x,y \in X$. But $T$ does not have a fixed point. This leads us to ask the following question:\vskip 0.3 cm
 \noindent What the topological properties that must be satisfied in order to assure the existence of fixed points for nonexpansive mappings?
\vskip 0.3 cm
\noindent The studies began in the framework of hyperconvex metric spaces which have the binary ball intersection property (any family of closed balls such that every two of them intersect, must have nonempty intersection) and the existence of fixed points for nonexpansive mappings in this setting is assured. In the proof of the existence of this fixed points, the behavior of the collection of admissible subsets of a hyperconvex metric space will play a crucial role. For a good reading on this subject, we can quote, for example (Chapter 4 in \cite{Kha}). If we pass to the framework of Banach spaces, the importance of their geometry in the study of the existence of fixed points for nonexpansive mappings or for generalized nonexpansive mappings, was still unknown until the mid 60s with the appearance of the works by F. E. Browder, D. G\"ohde and W. A. Kirk (\cite{Bro1, Goh, Kir1}). In fact, the first two authors showed the following result:

\vskip 0.3 cm
\noindent If $K$ is a bounded closed convex subset of uniformly convex Banach space X and if $T:X \rightarrow X$ is a nonexpansive mapping then $T$ admits a fixed point in $K$. W. A. Kirk (see \cite{Kir1}) has extended the above result to the case of reflexive Banach spaces that have normal structure.

\vskip 0.3 cm
\noindent A Banach space $X$ is said to have normal structure (resp. weak normal structure) if every bounded closed (resp. weakly compact) convex subset $K$ of $X$ has normal structure (see \cite{Goe, Kha}). The weak* normal structure for dual Banach spaces can be defined by a same way when we consider weak* compact convex subsets (see for example \cite{Lim}). Notice that in the case of reflexive Banach spaces, the weak normal structure is nothing else but the normal structure and the weak* normal structure. This claim is not true in general in the case of nonreflexive Banach spaces. For example, the space $l_1$ has the weak normal structure but fails to posses normal structure. Also, it is important to recall that weak normal structure implies the w.f.p.p (weak fixed point property) associated to nonexpansive mappings but the converse is not true in general. Indeed, the space $c_0$ (the space of real sequences which converge to 0) has w.f.p.p but it does not have the weak normal structure (see \cite{Mau}). One of the important problems in this area which is still open until today is whether every reflexive Banach space has the w.f.p.p?

\vskip 0.3 cm
\noindent Several authors have studied the reciprocal of known fixed point theorems, that is to check whether the fixed point property of certain classes of mappings characterizes a topological property or a metric property. This question goes back to the work of C. Bassaga \cite{Bes} who showed that the fixed point property for contraction mappings on a metric space does not characterize the completeness. P. V. Subrahmanyam \cite{Sub} gave a positive answer for the characterization of completeness of a metric space by the fixed point property of Kannan mappings. In \cite{Lin}, P. K. Lin and Y. Sternfeld showed that for the case of normed spaces and for convex subsets, the compactness is characterized by the fixed point property of Lipschitzian mappings. In \cite{Dom}, T. Dominguez Benavides et al. showed that in the case of Banach spaces, for the bounded closed convex subsets, the weak compactness is characterized by the fixed point property of weakly continuous mappings which is considered to be the converse of Schauder-Tychonoff theorem. In \cite{Won}, C. S. Wong showed that the weak quasi-normal structure is characterized by the fixed point property of Kannan mappings. The example of the space $c_0$ showed in fact that the fixed point property of nonexpansive mappings does not necessarily characterize the weak normal structure. \vskip 0.3 cm
\noindent We start our study at the source, that is the metric framework. In section 2, we establish some fixed point results concerning mappings which diminish the radius of invariant admissible subsets and this is done in the framework of bounded hyperconvex metric spaces. We obtain the same result if we consider a bounded metric space $M$ with $\mathcal{A}(M)$ (the collection of admissible subsets) compact with normal structure. These results extend naturally the nonexpansive framework.
\vskip 0.3 cm
\noindent In section 3, we are interested with the Banach space framework. Here, the results of chapter 2 are applied where the convexity will play the same role as the admissibility property. The results obtained here extend those obtained by A. Amini-Harandi et al \cite{Ami} concerning wrt orbitally nonexpansive mappings and we show particularly that the fixed point property for mappings that diminish the radius of invariant convex subsets characterizes the weak normal structure. It is well known that the fixed point property of Kannan mappings characterizes the weak quasi-normal structure but it does not characterize the weak normal structure. In section 4, we show that the fixed point property of orbitally Kannan mappings that diminish the radius of orbits characterizes the weak normal structure and by means of this class, we give the necessary conditions under which a Banach space possesses  weak quasi-normal structure but not weak normal structure.
 \vskip 0.3 cm
\noindent In the end, we give final remarks especially on the main results of W. A. Kirk and N. Shahzad \cite{Kir} recently established, and on their validity if we consider the class of mappings which diminish the radius of invariant weakly closed subsets of a weakly compact (not necessarily convex) of a Banach space $X$.

\section{Admissible Subsets in Metric Spaces and Fixed Point Property}
\noindent We start this section by introducing the following notions and notations
\vskip 0.3 cm
\begin{definition}\label{def 1.1}\rm A metric space  $(M,d)$ is called hyperconvex if for any class of closed balls $\{B(x_{i},\delta_{i})\}_{i \in I}$ of $M$ satisfying that
\vskip 0.3 cm
\centerline{$d(x_{i}, x_{j})\leq r_{i}+ r_{j} \hskip 0.5 cm \ i,j \in I.$}
\vskip 0.3 cm
\noindent Then necessarily, we have $\displaystyle \bigcap_{i \in I}B(x_{i}, r_{i})\neq \emptyset.$
\end{definition}

\noindent For a subset $K$ of $(M,d)$, set
\begin{align}
r_{x}(K)= & \sup\{d(x,y): y\in K\} \nonumber \\
r_{M}(K)= & \inf\{r_{x}(K): x \in M\} \nonumber \\
r(K)=& \inf\{r_{x}(K): x \in K\} \nonumber \\
\delta(K)= & \sup\{d(x,y): (x,y)\in K^{2}\} \nonumber \\
C_{M}(K)= &  \{x \in M: r_{x}(K)= r(K)\} \nonumber \\
C(K)= & \{x \in K: r_{x}(K)= r(K)\} \nonumber \\
cov(K)= & \bigcap\{H_{i}, H_{i} \ \hbox{ is \ a \ closed \ ball \ and } \ K \subseteq H_i\} \nonumber
\end{align}

\noindent $r_{M}(K)$ is called the radius of $K$ (relative to $M$), $\delta(K)$ is the diameter of $K$, $C_{M}(K)$ is the center of $K$ (relative to $M$), $r(K)$ is the the chebyshev radius of $K$, $C(K)$ is called the Chebyshev center of $K$ and $cov(K)$ is called the cover of $K$.

\begin{remark}\label{rem 1.1}\rm It is well known that every metric hyperconvex space is complete (see Proposition 4.1 in \cite{Kha}).
\end{remark}

\noindent For a bounded subset $C$ of a hyperconvex metric space $M$, we can state the following fundamental properties:
\begin{enumerate}
  \item $cov(C)= \displaystyle \bigcap_{x \in M}B(x, r_{x}(C))$.
  \item $r_{x}(Cov(C))= r_{x}(C)$ for all $x \in M$.
  \item $r_{M}(Cov(C))= r_{M}(C)$.
  \item $r_{M}(C)= \displaystyle \frac{1}{2}\delta(C)$.
  \item $\delta(cov(C))= \delta(C)$.
\end{enumerate}

\noindent Now, we give the definition of admissible metric spaces which will play an important role in the sequel.

\begin{definition}\label{def 1.2}\rm A bounded subset $K$ of a metric space is called admissible if $K= cov(K)$.
\end{definition}

\begin{remark}\label{rem 1.2}\rm The collection of all admissible subsets of a metric space $M$ is denoted by $\mathcal{A}(M)$. A set $C_{0} \subset M$ is admissible if and only if $C_{0}$ can be written as the intersection of closed balls in $M$. As a consequence, we have $\mathcal{A}(M)$ is closed by arbitrary intersections. This property implies that each admissible set of a hyperconvex metric space is also hyperconvex.
\end{remark}

\noindent Now, we are in a position to give the definition of mappings which diminishes the radius of admissible sets of a metric space.

\begin{definition}\label{def 1.3}\rm Let $M$ be a metric space and let $T: M \longrightarrow M$ be a mapping. $T$ is said to be diminishing the radius of invariant admissible subsets of $M$ if for any $A \in \mathcal{A}(M)$ such that $A$ is $T$-invariant, we have
\vskip 0.3 cm
\centerline{$r_{T_{x}}(T(A))\leq r_{x}(A)$ for every $x \in M$.}
\vskip 0.3 cm
\end{definition}

\begin{definition}\label{def 1.4}\rm Let $T$ be a self-mapping on a metric space $M$. $T$ is called nonexpansive if
\vskip 0.3 cm
\centerline{$d(Tx,Ty)\leq d(x,y) \ $for all$ \ x,y \in M.$}
\vskip 0.3 cm
\end{definition}

\begin{definition}\label{def 1.5}\rm Let $T$ be a self-mapping on a metric space $M$. $T$ is called orbitally wrt nonexpansive if
\vskip 0.3 cm
\centerline{$d(Tx,Ty)\leq r_{x}(O_{T}(y)) \ $for all$ \ x,y \in M,$}
\vskip 0.3 cm
\noindent where $O_{T}(y)$ is the orbit of $y$ relative to $T$ which is defined by
\vskip 0.3 cm
\centerline{$O_{T}(y)= \{y, Ty, T^{2}y,.....\}.$}
\end{definition}

\begin{remark}\label{rem 1.3}\rm It is easy to observe that every nonexpansive mapping is orbitally wrt nonexpansive while the converse is in general not true as the following example shows:
\end{remark}

\begin{example}\label{exa 1.1}\rm (see \cite{Nic}) Let $M= [0,1]$ equipped with the usual norm $|.|$. Then $T: M \longrightarrow M$ defined by
\vskip 0.3 cm
\centerline{$Tx= \left\{
               \begin{array}{ll}
                 \displaystyle \frac{x}{2} & \hbox{if} \ x \in [\frac{1}{2},1];\\
                 \displaystyle \frac{x}{4} & \hbox{if}  \ x \in [0, \frac{1}{2}[.
               \end{array}
             \right.
$}
\vskip 0.3 cm
\noindent is orbitally wrt nonexpansive but it is not nonexpansive (since it is not continuous at the point $x_{0}= \displaystyle \frac{1}{2}$). So, one of the constraints of orbitally wrt nonexpansive mappings is that they are not necessarily continuous.
\end{example}

\begin{lemma}\label{lem 1.1}\rm Let $M$ be a metric space and let $T: M \longrightarrow M$ be an orbitally wrt nonexpansive mapping. Then $T$ diminishes the radius of invariant admissible subsets of $M$.
\end{lemma}

\proof Denote by $\mathcal{H}$ the family
\vskip 0.3 cm
\centerline{$\mathcal{H}= \{K \in \mathcal{A}(M), K \neq \emptyset \ $and$ \ K \ \hbox{is} \ $T$-\hbox{invariant}\}.$}
\vskip 0.3 cm
\noindent It is easy to observe that $\mathcal{H}\neq \emptyset$ (since $M \in \mathcal{H}$). Let $K \in \mathcal{H}$, since $T$ is orbitally wrt nonexpansive, we have
\vskip 0.3 cm
\centerline{$d(Tx, Ty)\leq r_{x}(O_{T}(y)) \ $for all$ \ x \in M \ $and all$ \ y \in K.$}
\vskip 0.3 cm
\noindent On the other hand, the fact that $O_{T}(y) \subset K$ implies
\vskip 0.3 cm
\centerline{$r_{x}(O_{T}(y))\leq r_{x}(K).$}
\vskip 0.3 cm
\noindent Consequently, by passing to the $\sup$, it follows that
\vskip 0.3 cm
\centerline{$r_{T_{x}}(T(K))= \displaystyle \sup_{y \in K} d(Tx,Ty)\leq r_{x}(K),$}
\vskip 0.3 cm
\noindent which is the desired result.

\begin{corollary}\label{cor 1.1}\rm Let $M$ be a metric space and let $T: M \longrightarrow M$ be a nonexpansive mapping. Then $T$ diminishes the radius of invariant admissible subsets of $M$.
\end{corollary}

\noindent For a self-mapping $T$ on a metric space $M$, $F(T)$ will denote the set of fixed points of $T$, in other words;
\vskip 0.3 cm
\centerline{$F(T)= \{x \in M / Tx= x\}.$}
\vskip 0.3 cm
\noindent Our first fixed point result in this section is given by the following

\begin{theorem}\label{th 1.1}\rm Let $(M,d)$ be a bounded hyperconvex metric space and let $T: M \longrightarrow M$ be a mapping which diminishes the radius of invariant admissible subsets of $M$. Then $F(T)\neq \emptyset$, in addition, $F(T)$ is hyperconvex.
\end{theorem}

\proof The proof of the first claim can be obtained by an adaptation of the techniques of Theorem 4.8 in \cite{Kha}. Now, we will prove that $F(T)$ is hyperconvex. Let $\{B(x_{i}, \delta_{i})\}_{i \in I}$ be a family of closed balls centered at points $x_{i} \in F(T)$ such that
\vskip 0.3 cm
\centerline{$d(x_{i},x_{j}) \leq \delta_{i}+ \delta_{j}.$}
\vskip 0.3 cm
\noindent Since $M$ is hyperconvex, we get
\vskip 0.3 cm
\centerline{$K= \displaystyle \bigcap_{i \in I}B(x_{i}, \delta_{i})\neq \emptyset.$}
\vskip 0.3 cm
\noindent Here, we have the following cases:
\begin{itemize}
  \item $F(T)= \{x_{0}\} (x_{0}\in M)$. In this case, we have $\displaystyle \bigcap_{i \in I}B(x_{0}, \delta_{i})\bigcap F(T)= \{x_{0}\}\neq \emptyset$.
  \item If there exists $k_{0} \in I$ such that $B(x_{k_{0}}, \delta_{k_{0}})= \{x_{k_{0}}\}$, then
  \vskip 0.3 cm
  \centerline{$K \bigcap F(T)= \{x_{k_{0}}\}.$}
  \vskip 0.3 cm
  \noindent Therefore $F(T)$ is hyperconvex.
  \item Assume that $\forall i \in I, B(x_{i}, \delta_{i})$ and $F(T)$ have more than one point. Then each $B(x_{i}, \delta_{i})$ is $T$-invariant. Indeed, by taking $A = \{x_i\}$ and $z$ an arbitrary element in $B(x_i, \delta_i)$, we get
  \vskip 0.3 cm
  \centerline{$r_{Tz}(T\{x_i\}) = \|Tz - Tx_i\| = \|Tz - x_i\| \leq r_{z}(\{x_i\}) = \|z - x_i\| \leq \delta_{i}$}
  \vskip 0.3 cm
  which implies that $T(B(x_{i}, \delta_{i}))\subseteq B(x_{i}, \delta_{i})$. Consequently, we obtain that
  \vskip 0.3 cm
  \centerline{$T(K)= T(\displaystyle \bigcap_{i \in I}B(x_{i}, \delta_{i}))\subseteq \displaystyle \bigcap_{i \in I}T(B(x_{i}, \delta_{i}))\subseteq \displaystyle \bigcap_{i \in I}B(x_{i}, \delta_{i})= K.$}
  \vskip 0.3 cm
  This proves that $K$ is $T$-invariant.
  \item If $K= \{z_{0}\} (z_{0} \in M)$ then $z_{0}\in F(T)$ and hence $K \bigcap F(T)= \{z_{0}\}$. So, our claim holds.
  \item If $\delta(K)> 0$, then since $K$ is admissible and $T$-invariant and the fact that each admissible subset of a hyperconvex space is hyperconvex (see Remark 2.2), we deduce that $K$ is hyperconvex. By using the first assertion of this theorem, it follows that
      \vskip 0.3 cm
      \centerline{$K \bigcap F(T) \neq \emptyset$}
      \vskip 0.3 cm
      which is the desired result.
\end{itemize}

\noindent As an immediate consequence of Theorem \ref{th 1.1}, we have

\begin{corollary}\label{cor 1.2}\rm Let $(M,d)$ be a bounded hyperconvex metric space and let $T: M \longrightarrow M$ be an orbitally wrt nonexpansive mapping. Then $F(T)\neq \emptyset$ and $F(T)$ is hyperconvex.
\end{corollary}

\begin{corollary}\label{cor 1.3}\rm Let $(M,d)$ be a bounded hyperconvex metric space and let $T: M \longrightarrow M$ be a nonexpansive mapping. Then $F(T)\neq \emptyset$ and $F(T)$ is hyperconvex.
\end{corollary}

\begin{definition}\label{def 1.6}\rm Let $(M,d)$ be a metric space. The family of admissible subsets $\mathcal{A}(M)$ is said to be compact if every descending chain of nonempty elements of $\mathcal{A}(M)$ has a nonempty intersection.
\end{definition}

\begin{definition}\label{def 1.7}\rm $\mathcal{A}(M)$ is said to have normal structure if $r(K) < \delta(K)$ for every $K \in \mathcal{A}(M)$ with $\delta(K)> 0$
\end{definition}

\begin{remark}\label{rem 1.4}\rm It can be easily seen that $\mathcal{A}(M)$ have normal structure if and only if for every $K \in \mathcal{A}(M)$ with $\delta(K)> 0$ there exists $r_{0}> 0$ and $z_{0} \in K$ such that
\vskip 0.3 cm
\centerline{$K \subseteq B(z_{0}, r_{0}).$}
\vskip 0.3 cm
\end{remark}

\noindent Our main result in this section is given by the following

\begin{theorem}\label{th 1.2}\rm Let $(M,d)$ be a bounded metric space such that $\mathcal{A}(M)$ is compact and having normal structure and let $T: M \longrightarrow M$ be a mapping which diminishes the radius of invariant admissible subsets of $M$, then $T$ has a fixed point in $M$.
\end{theorem}

\proof Following the compactness of $\mathcal{A}(M)$ and using Zorn's lemma, we can find a subsets $K_{0} \in\mathcal{A}(M)$ such that $K_{0}$ is minimal (relative to the inclusion) which is $T$-invariant. It is easy to observe that
\vskip 0.3 cm
\centerline{$T(K_{0})\subseteq \displaystyle \bigcap_{x \in M}B(x, r_{x}(T(K_{0})))$}
\vskip 0.3 cm
\noindent and the fact that $r_{x}(T(K_{0}))\leq r_{x}(K_{0})$ (since $T(K_{0})\subset K_{0}$), gives
\vskip 0.3 cm
\centerline{$cov(T(K_{0}))\subset \displaystyle \bigcap_{x \in M}B(x, r_{x}(T(K_{0})))\subset \displaystyle \bigcap_{x \in M}B(x, r_{x}(K_{0}))= K_{0}.$}
\vskip 0.3 cm
\noindent Then
 \vskip 0.3 cm
\centerline{$T(cov(T(K_{0})))\subset T(K_0)\subseteq cov(T(K_{0})),$}
\vskip 0.3 cm
\noindent which implies that $cov(T(K_{0}))$ is $T$-invariant. The minimality of $K_{0}$ implies that
\vskip 0.3 cm
\centerline{$K_{0}= cov (T(K_{0})).$}
\vskip 0.3 cm
\noindent Since $\mathcal{A}(M)$ has normal structure, it follows that
\vskip 0.3 cm
\centerline{$\widetilde{K_{0}}= \{x \in K_{0}: K_{0}\subseteq B(x,r_{0})\}\neq \emptyset,$}
\vskip 0.3 cm
\noindent where
\vskip 0.3 cm
\centerline{$r(K_{0})< r_{0}< \delta(K_{0}).$}
\vskip 0.3 cm
\noindent By a simple observation, we get
\vskip 0.3 cm
\centerline{$\widetilde{K_{0}}= \displaystyle \bigcap_{x \in K_0}B(x,r_{0})\bigcap K_{0}$}
\vskip 0.3 cm
\noindent which prove that
\vskip 0.3 cm
\centerline{$\widetilde{K_{0}}\in \mathcal{A}(M).$}
\vskip 0.3 cm
\noindent Next, for all $z_{0}\in \widetilde{K_{0}}$ and all $x \in K_{0}$, we have
\vskip 0.3 cm
\centerline{$d(Tz_0,Tx)\leq r_{Tz_0}(T(K_{0}))\leq r_{z_0} (K_{0})= \displaystyle \sup_{y \in K_{0}}d(z_0,y)\leq r_{0}.$}
\vskip 0.3 cm
\noindent Hence $Tx \in B(Tz_{0}, r_{0})$ for every $x \in K_{0}$ which implies that
\vskip 0.3 cm
\centerline{$T(K_{0})\subseteq B(Tz_{0}, r_{0}).$}
\vskip 0.3 cm
\noindent Consequently,
\vskip 0.3 cm
\centerline{$K_{0}= cov(T(K_{0}))\subseteq B(T(z_{0}), r_{0}).$}
\vskip 0.3 cm
\noindent This shows that
\vskip 0.3 cm
\centerline{$Tz_{0}\in \widetilde{K_{0}}.$}
\vskip 0.3 cm
\noindent In other words, the set $\widetilde{K_{0}}$ is $T$-invariant. But
\vskip 0.3 cm
\centerline{$\delta(\widetilde{K_{0}})\leq r < \delta(K_{0}).$}
\vskip 0.3 cm
\noindent Hence $\widetilde{K_{0}} \varsubsetneq K_{0}$. In addition, since $\widetilde{K_{0}}$ is $T$-invariant with $\widetilde{K_{0}}\in \mathcal{A}(M)$, we obtain a contradiction with the minimality of $K_{0}$. Hence, necessarily $\delta(K_{0})=0$ and $K_{0}$ consists of a point which is a fixed point of $T$.
\vskip 0.3 cm
\noindent The following corollaries follow
\begin{corollary}\label{cor 1.3}\rm Let $(M,d)$ be a bounded metric space such that $\mathcal{A}(M)$ is compact with normal structure. Then every orbitally wrt nonexpansive mapping $T: M \longrightarrow M$ has a fixed point in $M$.
\end{corollary}

\begin{corollary}\label{cor 1.4}\rm Let $(M,d)$ be a bounded metric space such that $\mathcal{A}(M)$ is compact and having normal structure. Then every nonexpansive mapping $T: M \longrightarrow M$ has a fixed point in $M$.
\end{corollary}

\begin{definition}\label{def 1.8}\rm The family $\mathcal{A}(M)$ is said to have quasi-normal structure if for each $K \in \mathcal{A}(M)$ with $\delta(K)> 0$ there exists $x_{0}\in K$ such that for all $x \in K$,
\begin{equation}
d(x_{0},x )< \delta(K).
\end{equation}
\end{definition}

\begin{remark}\label{rem 1.5}\rm It is easy to observe that normal structure implies quasi-normal structure and the converse is not true in general. The interest of quasi-normal structure is well understood in the setting of Banach spaces as we will show in the next section.
\end{remark}

\begin{definition}\label{def 1.9}\rm A mapping $T: M \longrightarrow M$ is called a Kannan mapping if
\begin{equation}
d(Tx,Ty)\leq \displaystyle \frac{1}{2}(d(x,Tx)+ d(y,Ty)) \hskip 0.5 cm \hbox{for \ all} \  x,y \in M.
\end{equation}
\end{definition}

\noindent The link between quasi-normal structure and the fixed point property for Kannan mappings is given by the following theorem.
\begin{theorem}\label{th 1.3}\rm Assume that $(M,d)$ is a bounded metric space. If $\mathcal{A}(M)$ is countably compact having quasi-normal structure then every Kannan self-mapping on $M$ has a (unique) fixed point.
\end{theorem}

\noindent Conversely, if $\mathcal{A}(M)$ does not have quasi-normal structure then there exists $K_{0}\in \mathcal{A}(M)$ and a Kannan mapping $T: K_{0} \longrightarrow K_{0}$ which is fixed point free.

\section{Normal structure and mapping which diminish the radius of n.c.i (nontrivial convex invariant) subsets}
\noindent The aim of this section is to apply the results of the previous section to obtain fixed point results for self-mappings which diminish the radius of invariant convex subsets of weakly (resp. weak$^{\star}$) compact convex subsets of Banach spaces (resp. dual Banach spaces).
\vskip 0.3 cm
\noindent In 1948, M. S. Brodskii and D. P. Milman \cite{Bro} have introduced the concept of normal structure.
\begin{definition}\label{def 2.1}\rm A convex subset $K$ of a Banach space (resp. dual Banach space) $X$ is said to have normal structure (resp. weak normal structure, weak$^{\star}$ normal structure) if for any bounded closed (resp. weakly closed, weak$^{\star}$-closed) convex subset $K_{0}$ of $K$ such that $\delta(K_{0}) > 0$, there exists $x_{0}\in K_{0}$ for which
\begin{equation}
r_{x_{0}}(K_{0})< \delta(K_{0})
\end{equation}
\noindent $x_{0}$ is called a nondiametral point of $K_{0}$.
\end{definition}

\begin{definition}\label{def 2.2}\rm A Banach space (resp. dual Banach space) $X$ is said to have weak normal structure (resp. weak$^{\star}$) if every nontrivial weakly (resp. weak$^{\star}$) compact convex subset of $X$ has normal structure.
\end{definition}

\noindent For a bounded convex subset $A$ of a Banach space $X$, we will adopt the same definition for $r_{x}(K), r(K)$ and $C(K)$, just we replace the distance by the norm.

\begin{remark}\label{rem 2.1}\rm If $K$ is nontrivial weakly (resp. weak$^{\star}$) compact convex subset of a Banach space $X$ (resp. dual space) having weak normal structure (resp. weak$^{\star}$) normal structure then $C(K)$ is a nonempty admissible subset of $K$ such that $C(K) \varsubsetneq K$.
\end{remark}

\noindent The following result characterizes the normal structure of bounded convex subsets in Banach spaces.
\begin{theorem}\label{th 2.1}\rm A nonempty bounded convex subset $K$ of Banach space $X$ has normal structure if and only if $K$ does not contain a non constant sequence $(x_{n})_{n}$ such that
\begin{equation}\label{eq 2.2}
\displaystyle \lim_{n \longrightarrow +\infty}d(x_{n+1}, co(\{x_{i}\}_{i=1}^{n}))= \delta(\{x_{i}\}_{i=1}^{\infty})
\end{equation}
\noindent where $co$ denotes the convex hull and $d$ is the distance induced by the norm.
\end{theorem}

\noindent A sequence $(x_{n})$ which satisfies $(3.2)$ is called a diametral sequence.

\begin{definition}\label{def 2.3}\rm A Banach space (resp. dual Banach space) $X$ is said to have weak fixed point property (in short; w.f.p.p) (resp. weak$^{\star}$ fixed point property) (in short; w$^{\star}$.f.p.p) if for every nontrivial weakly (resp. weak$^{\star}$) compact convex subset $K$ of $X$, every nonexpansive mapping $T: K \longrightarrow K$ has at least a fixed point in $K$.
\end{definition}

\begin{remark}\label{rem 2.2}\rm It is well known that weak (resp. weak$^{\star}$)  normal strcuture implies w.f.p.p (resp. w$^{\star}$.f.p.p) while the converse is not true in general as the following example shows:
\end{remark}

\begin{example}\label{exa 2.1}\rm (see page 1 in \cite{Kar}) Let $X= l^{2}_{|.|}$ the Banach space $l^{2}$ equipped with the norm $|.|$ defined by
\vskip 0.3 cm
\centerline{$|x|= \max\{\|x\|_{\infty}, \displaystyle \frac{\|x\|_{2}}{\sqrt{2}}\}$}
\vskip 0.3 cm
\noindent where $\|.\|_{\infty}$ is the $l^{\infty}$ norm and $\|.\|_{2}$ is the euclidean norm in $l^{2}$. It was proved that $X$ has w.f.p.p. However, $X$ does not have weak normal structure since all the points of the bounded convex subset
\vskip 0.3 cm
\centerline{$K_{0}= \{ x = (x_{n})_{n} : x_{n}\geq 0 \ $for all$ \ n \in \mathbb{N}, \|x\|_{2}\leq 1\}$}
\vskip 0.3 cm
\noindent are diametral.
\end{example}

\noindent \textbf{\it Question}: For an arbitrary Banach space $X$, is it possible to find a class of mappings which contains "strictly" nonexpansive mappings such that the fixed point property associated to this class characterizes the weak normal structure?

\vskip 0.3 cm
\noindent For a nontrivial weakly (resp weak$^{\star}$) compact convex subset $K$ of a Banach space (resp. dual Banach space) $X$, we denote by $\mathcal{F}_{K}$ the set of self-mappings on $K$ which diminish the radius of invariant convex subsets of $K$.

\vskip 0.3 cm
\begin{definition}\label{def 2.5}\rm A Banach space (resp dual Banach space) $X$ is said to have the weak (resp. weak$^{\star}$) fixed point property for mappings which diminish the radius of invariant convex subsets if for every nontrivial weakly (resp weak$^{\star}$) compact convex subset $K$ of $X$, every $T \in \mathcal{F}_{K}$ has at least a fixed point in $K$.
\end{definition}

\noindent Our main result in this section is given by the following theorem

\begin{theorem}\label{th 2.2}\rm Let $X$ be a Banach space. Then $X$ has weak normal structure if and only if $X$ has the fixed point property for mappings which diminish the radius of invariant convex subsets.
\end{theorem}

\proof \textbf{\it Necessary condition}
\vskip 0.3 cm
\noindent It is an immediate consequence of Theorem 2.2 and Remark 3.1.
\vskip 0.3 cm
\noindent \textbf{\bf \it Sufficient condition}
\vskip 0.3 cm
\noindent Assume that $X$ has not weak normal structure then there exists a nontrivial weakly compact convex subset $K_{0}$ of $X$ such that $K_{0}$ does not have weak normal structure. Following Theorem 3.1, $K_{0}$ contains a diametral nonconstant sequence $\{x_{n}\}$ such that
\begin{equation}\label{eq 2.3}
\displaystyle \lim_{n \longrightarrow + \infty}d(x_{n+1}, co(\{x_{n}\}_{i =1}^{n}))= \delta(\{x_{n}\}_{n})> 0.
\end{equation}
\noindent Denote by $\widetilde{K_{0}}= \overline{co}(\{x_{n}\}_n)$ where $\overline{co}$ is the closed convex hull. Let $T: \widetilde{K_{0}} \longrightarrow \widetilde{K_{0}}$ defined as in \cite{Ami} by
\vskip 0.3 cm
\centerline{$Tx= \left\{
                 \begin{array}{ll}
                   x_{1} & \hbox{if} \ x \notin \{x_{n}, n \in \mathbb{N}\}; \\
                   x_{n+1} & \hbox{if} \ x = x_{n} \ \hbox{for \ some} \ n.
                 \end{array}
               \right.
$}
\vskip 0.3 cm
\noindent Now, let $A$ be convex subset of $\widetilde{K_{0}}$ which is $T$-invariant. Then, we have for all $x \in \widetilde{K_0}$ and all $y \in \widetilde{K_{0}}$,
\vskip 0.3 cm
\centerline{$r_{T_{x}}(T(A))= r_{y}(A)= \delta(\widetilde{K_{0}})> 0$.}
\vskip 0.3 cm
\noindent Consequently,  by taking $y = x$, we deduce that $T$ diminishes the radius of invariant convex subsets. However, $T$ is a fixed point free mapping.

\begin{remark}\label{rem 2.3}\rm In the previous example, if $x \notin \{x_{n}, n \in \mathbb{N}\}$ then $co(O_{T}(x))$ is a good example of a nontrivial convex $T$-invariant subset of $\widetilde{K_{0}}$. Indeed, let $z \in co(O_{T}(x))$.  If $z \notin \{x_{n}, n \in \mathbb{N}\}$ then $Tz= x_{1}\in O_{T}(x)\subseteq co(O_{T}(x))$. If $z=x_{n_0}$ for some $n_{0}$ then we have $Tx_{n_{0}}= x_{n_{0}+ 1}\in O_{T}(x)\subseteq co(O_{T}(x))$.
\end{remark}

\noindent Since every uniformly convex Banach space has the weak normal structure, we can derive the following

\begin{corollary}\label{cor 2.1}\rm Every uniformly convex Banach space has the weak fixed point property for mappings which diminish the radius of invariant convex subsets.
\end{corollary}

\begin{corollary}\label{cor 2.2}\rm If $X$ is one the following Banach spaces
\begin{enumerate}
 \item $l^{1}$ space;
  \item $(J, \|.\|_{_{j}}); \ j=1,2, $ (see \cite{Kuc})
 \item James tree space $JT$. (see \cite{Kuc})
\end{enumerate}
\noindent Then $X$ has the weak$^{\star}$ fixed point property for mappings which diminish the radius of invariant convex subsets.
\end{corollary}

\proof This result follows from the fact that these spaces has weak$^{\star}$ normal structure.
\vskip 0.3 cm
\noindent Since weakly compact subsets of dual spaces are weak$^{\star}$ compact convex then we can deduce the following

\begin{corollary}\label{cor 2.3}\rm If $X$ is one of the following Banach spaces
\begin{enumerate}
  \item $l^{1}$ space;
  \item $(J, \|.\|_{_{j}}), j=1,2 $;
  \item James tree space JT.
\end{enumerate}
\noindent Then $X$ has the weak fixed point property for mappings which diminish the radius of invariant convex subsets.
\end{corollary}

\vskip 0.3 cm

\section{Normal Structure and Orbitally Kannan Mappings}
\noindent First of all, let us give the following definitions

\begin{definition}\label{def 3.1}\rm Let $C$ be a nonempty subset of a Banach space $X$ and let $T: C \longrightarrow C$ be a mapping. $T$ is said to be orbitally Kannan mapping if
\begin{equation}\label{eq 3.1}
\|Tx -Ty\|\leq \displaystyle \frac{1}{2}(r_{x}(O_{T}(x))+ r_{y}(O_{T}(y))) \hbox{ for \ all} \ x,y \in C.
\end{equation}
\end{definition}

\begin{definition}\label{def 3.2}\rm Let $C$ be a nonempty subset of a Banach space $X$ and let $T: C \longrightarrow C$ be a mapping. $T$ is said to be diminishing the radius of orbits if for all $x \in C$, we have
\begin{equation}\label{eq 3.2}
r_{T{x}}(O_{T}(Tx))\leq  r_{x}(O_{T}(x)).
\end{equation}
\noindent Kannan mappings are given by Definition 2.9 where the distance is replaced by the norm.
\end{definition}

\begin{lemma}\label{lem 3.1}\rm Let $C$ be a nonempty subset of a Banach space $X$ and let $T: C \longrightarrow C$ be a Kannan self-mapping. Then $T$ is orbitally Kannan which diminishes the radius of orbits.
\end{lemma}

\proof Let $T: C \longrightarrow C$ be a Kannan mapping then
\begin{equation}
\|Tx-Ty\|\leq \displaystyle \frac{1}{2}(\|x- Tx\|+ \|y-Ty\|) \ \hbox{for \ all} \ x,y \in C.
\end{equation}
\noindent Now, since
\begin{equation}\label{eq 3.4}
\|z- Tz\|\leq r_{z}(O_{T}(z)) \ \hbox{for \ all} \ z \in C,
\end{equation}
\noindent we get
\begin{equation}\label{eq 3.5}
\|Tx -Ty\|\leq \displaystyle \frac{1}{2}(r_{x}(O_{T}(x))+ r_{y}(O_{T}(y))) \ \hbox{for \ all} \ x,y \in C,
\end{equation}
\noindent which proves that $T$ is orbittally Kannan.
\vskip 0.3 cm
\noindent Next, for all $x \in C$ and for all integer $n \geq 1$, we have
\begin{equation}\label{eq 3.6}
\|Tx -T^{n}x\|\leq \displaystyle\frac{1}{2} (\|x- Tx\|+ \|T^{n-1}x-T^{n}x\|).
\end{equation}
\noindent It is easy to prove that the sequence $(\|T^{n-1}x- T^{n}x\|)_{n}$ is decreasing, so for all $n \geq 1$
\begin{equation}\label{eq 3.7}
\|Tx -T^{n}x\|\leq \|x- Tx\| \leq r_{x}(O_{T}(x)).
\end{equation}
\noindent Using (4.5), we get
\begin{equation}\label{eq 3.8}
\|Tx -T^{n}x\|\leq \displaystyle\frac{1}{2} (r_{x}(O_{T}(x))+ r_{x}(O_{T}(x))),
\end{equation}
\noindent Consequently, by passing to the $\sup$ (relative ton $n$), it follows that
\begin{equation}\label{eq 3.9}
r_{Tx}(O_{T}(Tx))\leq r_{x}(O_{T}(x)).
\end{equation}
\noindent which completes the proof
\vskip 0.3 cm
\noindent The following example shows that, there exist orbitally Kannan mappings which diminish the radius of orbits which are not of Kannan type.

\begin{example}\label{exa 3.1}\rm Let
\begin{align}
T: [0,1] & \longrightarrow [0,1] \nonumber \\
& x \longrightarrow \displaystyle \frac{x}{2} \nonumber
\end{align}
\noindent $T$ is not a Kannan mapping since the inequality
\begin{align}
|Tx- Ty|=& |\displaystyle \frac{x}{2}- \displaystyle \frac{y}{2}| \leq \displaystyle\frac{1}{2}(|x- \displaystyle \frac{x}{2}|+ |y- \displaystyle \frac{y}{2}|)= \displaystyle \frac{1}{2}(\displaystyle \frac{x}{2}+ \displaystyle \frac{y}{2}) \nonumber \\
=& \displaystyle \frac{1}{4}(x+ y) \nonumber
\end{align}
\noindent is not satisfied for all $x,y \in [0,1]$. To see this, it suffices to take $x_{0}= 1$ and $y_{0}= 0$.
\vskip 0.3 cm
\noindent On the other hand, for all $x \in [0,1]$, we have $r_{x}(O_{T}(x))= \displaystyle \sup_{n \geq 1}|x- \displaystyle \frac{x}{2^{n}}|= x$, which gives that
\begin{equation}\label{eq 3.10}
|Tx- Ty|= \displaystyle \frac{1}{2}|x- y|\leq \displaystyle \frac{1}{2}(r_{x}(O_{T}(x))+ r_{y}(O_{T}(y))= \displaystyle \frac{1}{2}(x+ y)
\end{equation}
\noindent Furthermore, for all integer $n \geq 1$,
\begin{align}
|Tx- T^{n}x|= |\displaystyle \frac{x}{2}- \displaystyle \frac{x}{2^{n}} |= \displaystyle \frac{1}{2}|x- \displaystyle \frac{x}{2^{n-1}}|\leq & \displaystyle \frac{1}{2}\sup_{n \geq1}|x- \displaystyle \frac{x}{2^{n}}| \nonumber \\
\leq & \displaystyle \frac{1}{2} r_{x}(O_{T}(x)) \nonumber \\
\leq & r_{x}(O_{T}(x)) \nonumber
\end{align}
\end{example}

\noindent By passing to the $\sup$ (relative ton $n$), we get
\begin{equation}
r_{Tx}(O_{T}(Tx))\leq r_{x}(O_{T}(x))
\end{equation}
\noindent which completes the proof.

\begin{definition}\label{def 3.3}\rm Let $C$ be nontrivial weakly compact convex subset of a Banach space $X. \ C$ is said to be quasi-normal if for any nontrivial closed convex subset $K_{0}$ of $C$, there exists $x_{0} \in K_{0}$ such that
\begin{equation}
\|x_{0}- y\| < \delta(K_{0}), \hbox{for all} \ y \in K_{0}.
\end{equation}
\end{definition}

\begin{definition}\label{def 3.4}\rm A Banach space $X$ is said to have weak quasi-normal structure if every nontrivial weakly compact convex subset of $X$ has quasi-normal structure.
\end{definition}

\begin{remark}\label{def 3.1}\rm Similar to Remark 2.5, weak normal structure implies weak quasi-normal structure while the converse is not true in general.
\end{remark}

\noindent In \cite{Won} and \cite{Won1}, C. S Wong has established the following important results. The first one characterizes the weak quasi-normal structure for weakly compact convex subsets in Banach spaces.

\vskip 0.3 cm
\begin{theorem}\label{th 3.1}\rm Let $X$ be a Banach space. Then the following two assertions are equivalent
\begin{description}
  \item[$(\imath)$] Every weakly compact convex subset $X$ has quasi-normal structure;
  \item[$(\imath\imath)$] Every Kannan self-mapping on a weakly compact convex subset of $X$ has a (unique) fixed point.
\end{description}
\end{theorem}

\begin{theorem}\label{th 3.2}\rm If $X$ in one of the following Banach spaces
\begin{description}
  \item[$(\imath)$] $X$ separable;
  \item[$(\imath\imath)$] $X$ strictly convex;
  \item[$(\imath\imath)$] $X$ has Kadee-Klee property, that is for any sequence $(x_{n})$ in $X$, if $x_{n}$ converges weakly to some $x \in X$ and $\|x_{n}\| \longmapsto \|x\|$ then $(x_{n})_{n}$ converges in norm to $x$.
\end{description}
\noindent Then $X$ has weak quasi-normal structure.
\end{theorem}

\begin{remark}\label{rem 3.2}\rm Theorem 4.1 can be seen as the variant of Theorem 2.3 in the case of Banach spaces.
\end{remark}

\begin{remark}\label{rem 3.3}\rm By Theorem 4.2, the Lebesgue space $L^{1}([0,1])$ has the weak quasi-normal structure. However, this space fails to have the weak normal structure due to the famous Alspach's result (see \cite{Alp}).
\end{remark}

\begin{remark}\label{rem 3.4}\rm From Remark 4.3 and Theorem 4.2, it is easy to deduce that the w.f.p.p for Kannan mappings does not necessarily characterize weak normal structure.
\end{remark}

\noindent One of our goals is to show that weak normal structure is characterized by the w.f.p.p for orbitally Kannan mappings which diminish the orbits.
\vskip 0.3 cm
\noindent We start our investigation by the following result.
\begin{theorem}\label{th 3.3}\rm Let $C$ be a weakly compact convex subset of a Banach space $X$ and let $T: C \longrightarrow C$ be an orbitally Kannan mapping which diminishes the radius of orbits. Then the set
\vskip 0.3 cm
\centerline{$\{r_{x}(O_{T}(x)), x \in C\}$}
\vskip 0.3 cm
\noindent has a smallest number.
\end{theorem}

\proof Set
\vskip 0.3 cm
\centerline{$C_{\delta} =\{z \in C: r_{z}(O_{T}(z))\leq \delta\}, \hskip 0.3 cm \delta > 0.$}
\vskip 0.3 cm
\noindent Since $C$ is weakly compact then $C$ is bounded, hence there exists $\delta_{0}> 0$ such that the set $C_{\delta_0}$ is nonempty. Our aim now is to prove that
\vskip 0.3 cm
\centerline{$\displaystyle \bigcap_{\delta \in J}C_{\delta}\neq \emptyset$}
\vskip 0.3 cm
\noindent where
\vskip 0.3 cm
\centerline{$J= \{\delta > 0: C_{\delta}\neq \emptyset\}.$}
\vskip 0.3 cm
\noindent We denote by $\overline{co}(T(C_{\delta}))$ the closed convex hull of the set $T(C_{\delta})$.
\vskip 0.3 cm
\noindent Next, it can be easily seen that the family $\{C_{\delta}, \delta \in J\}$ is increasing relative to the usual order of $J$. The weak compactness of $C$ implies that the family $\{\overline{co}(T(C_{\delta})): \delta \in J\}$ has the finite intersection property. To complete the proof, it suffices to show that $\overline{co}(T(C_{\delta}))\subseteq C_{\delta}$ for all $\delta \in J$. Let $r_{0}$ an arbitrary element of $J$ and let $x \in \overline{co}(T(C_{r_{0}}))$, fix $\epsilon > 0$, then there exist $\alpha_1, ..., \alpha_n \in [0, 1]$ and $x_1, x_2,...,x_n \in C_{r_{0}}$ such that $\displaystyle \sum_{i=1}^{n}\alpha_{i}= 1$ and
\vskip 0.3 cm
\centerline{$\|x- \displaystyle \sum_{i=1}^{n}\alpha_{i} x_{i}\| < \displaystyle \frac{\epsilon}{2}$}
\vskip 0.3 cm
\noindent Then by the triangle inequality, for all integer $k_{0}\geq 1$, we have
\begin{align}
\|x- T^{k_{0}}x\|\leq & \|x- \displaystyle \sum_{i=1}^{n}\alpha_{i} T(x_{i})\|+ \|\displaystyle \sum_{i=1}^{n}\alpha_{i} T(x_{i})- T^{k_{0}}(x)\| \nonumber \\
< & \displaystyle \frac{\epsilon}{2} + \displaystyle \sum_{i=1}^{n} \alpha_i \|T(x_{i})- T^{k_{0}}(x)\| \nonumber \\
< & \displaystyle \frac{\epsilon}{2} + \displaystyle \sum_{i=1}^{n}\frac{\alpha_{i}}{2}(r_{x_{i}}(O_{T}(x_{i}))+ r_{T^{k_0 - 1}x}(O_{T}(T^{k_0 - 1}x))). \nonumber \\
\end{align}

\noindent Since $x_{i}\in C_{r_{0}}$ for all $i=1,...,n$, by definition, we have
\vskip 0.3 cm
\centerline{$r_{x_{i}}(O_{T}(x_{i}))\leq r_{0}.$}
\vskip 0.3 cm
\noindent Consequently,
\vskip 0.3 cm
\centerline{$\displaystyle \sum_{i=1}^{n}\frac{\alpha_{i}}{2}r_{x_{i}}(O_{T}(x_{i}))\leq \displaystyle \frac{r_{0}}{2}.$}
\vskip 0.3 cm
\noindent On the other hand, since $T$ diminishes the radius of each orbit, it follows that
\vskip 0.3 cm
\centerline{$r_{T^{k_0 - 1}x} (O_T(T^{k_{0}-1}(x)))\leq r_{x}(O_{T}(x)).$}
\vskip 0.3 cm
\noindent Hence, for all $k_{0}\geq 1$, we get
\vskip 0.3 cm
\centerline{$\|x- T^{k_{0}}(x)\|< \displaystyle \frac{\epsilon}{2}+ \displaystyle \frac{r_{0}}{2}+ \displaystyle \frac{r_{x}(O_{T}(x))}{2}.$}
\vskip 0.3 cm
\noindent Passing to the $\sup$ (relative to $k_0$), we infer that
\vskip 0.3 cm
\centerline{$r_{x}(O_{T}(x))= \displaystyle \sup_{k_{0}\geq1}\|x- T^{k_{0}}(x)\| \leq \displaystyle \frac{\epsilon}{2} + \displaystyle \frac{r_{0}}{2} + \displaystyle \frac{r_{x}(O_{T}(x))}{2},$}
\vskip 0.3 cm
\noindent which proves that
\vskip 0.3 cm
\centerline{$r_{x}(O_{T}(x)) \leq \epsilon+ r_{0}$}
\vskip 0.3 cm
\noindent Thus,
\vskip 0.3 cm
\centerline{$r_{x}(O_{T}(x)) \leq  r_{0}.$}
\vskip 0.3 cm
\noindent This implies that $x \in C_{r_{0}}$ which is the desired result.

\begin{theorem}\label{th 3.4}\rm Let $C$ be a weakly compact convex subset of a Banach space $X$. Let $T$ be an orbitally Kannan mapping on $C$ which diminishes the radius of each orbit. Then the following assertions are equivalent:
\begin{enumerate}
  \item $T$ has a unique fixed point;
  \item $\displaystyle \inf\{r_{x} (O_T(x)): x \in C\}= 0$;
  \item for all $x \in C$ with $r_{x} (O_T(x))> 0$, there exists $y \in C$ such that
\vskip 0.3 cm
\centerline{$r_{x}(O_{T}(x))> r_{y}(O_{T}(y)) \ $for some$ \ y \in C.$}
\end{enumerate}
\end{theorem}

\proof $1) \Longrightarrow 2)$
\vskip 0.3 cm
\noindent If $T$ has a unique fixed point $x_{0}\in C$ then $r_{x_{0}}(O_{T}(x_{0}))= 0$ and hence $\displaystyle \inf\{r_{x} (O_T(x)): x \in C\}= 0$.
\vskip 0.3 cm
\noindent $2) \Longrightarrow 3)$
\vskip 0.3 cm
\noindent Assume that there exists $z_{0}\in C$ such that $r_{z_{0}}(O_{T}(z_{0})) > 0$ and
\vskip 0.3 cm
\centerline{$r_{z_{0}}(O_{T}(z_{0}))\leq r_{y}(O_{T}(y)) \ $for all$ \ y \in C.$}
\vskip 0.3 cm
\noindent By passing to the $\inf$, it follows that
\vskip 0.3 cm
\centerline{$r_{z_{0}}(O_{T}(z_{0}))\leq \displaystyle \inf_{y \in C}\{r_{y}(O_{T}(y)): y \in C\}= 0.$}
\vskip 0.3 cm
\noindent Consequently, we obtain that
\vskip 0.3 cm
\centerline{$r_{z_{0}}(O_{T}(z_{0}))= 0,$}
\vskip 0.3 cm
\noindent which is a contradiction.
\vskip 0.3 cm
\noindent $3) \Longrightarrow 1)$
\vskip 0.3 cm
\noindent Assume that $T$ has no fixed points. Then for all $x \in C$, we have $r_{x}(O_{T}(x))> 0$. But following Theorem 4.3, there exists $y_{0}\in C$ such that $\inf\{r_{x}(O_{T}(x)): x \in C\}= r_{y_{0}}(O_{T}(y_{0})) > 0$. By assumption,
\vskip 0.3 cm
\centerline{$r_{y_{0}}(O_{T}(y_{0}))  > r_{x_{0}}(O_{T}(x_{0})) \ $for some$ \ x_{0}\in C, $}
\vskip 0.3 cm
\noindent which is a contradiction.

\begin{theorem}\label{th 3.5}\rm Let $C$ be a weakly compact convex subset of a Banach space $X$. Let $T: C \longrightarrow C$ be an orbitally Kannan mapping which diminish the radius of each orbit. Then
\begin{enumerate}
  \item Suppose that $T$ has no fixed points. Then there exists a $T$-invariant closed convex subset $K_{0}$ of $C$ such that $\delta(C_{0})> 0$ and $r_{x}(O_{T}(x))= \delta (C_{0})$ for all $x \in C_{0}$.
  \item If for each closed convex subset $K_0$ of $C$ with $\delta(K_0) > 0$ there exists $x_0 \in K_0$ such that
  \vskip 0.3 cm
  \centerline {$r_{x_0}(O_T(y)) < \delta(K_0) \ \hbox{for all} \ y \in K_0$.}
  \vskip 0.3 cm
  \noindent Then $T$ has a unique fixed point.
\end{enumerate}
\end{theorem}
\proof By Zorn's lemma, there exists a minimal $T$-invariant closed convex subset $\widetilde{K_0}$ of $C$. By Theorem 4.3, we assure the existence of an element $x_0 \in \widetilde{K_0}$ such that
\vskip 0.3 cm
\centerline {$\varepsilon_0 = r_{x_0}(O_T(x_0)) = \inf\{r_x (O_T(x)): x \in  \widetilde{K_0}\}.$}
\vskip 0.3 cm
\noindent Since $F(T) = \emptyset$, then by Theorem 4.4, we have $\varepsilon_0 > 0$. Using the same notations of Theorem 4.3, since $\overline{co}(T(C_{\epsilon_0})) \subset C_{\epsilon_0}$, it follows that

\vskip 0.3 cm
\centerline {$T(\overline{co}(T(C_{\epsilon_0}))) \subset T(C_{\epsilon_0}) \subset \overline{co} (T(C_{\epsilon_0})).$}
\vskip 0.3 cm
\noindent But $\widetilde{K_0}$ is minimal, so $\widetilde{K_0} = \overline{co} (T(C_{\epsilon_0}))$ and consequently,

\vskip 0.3 cm
\centerline {$C_{\epsilon_0} = \widetilde{K_0}.$}
\vskip 0.3 cm
\noindent By the choice of $\varepsilon_0$, we get
\vskip 0.3 cm
\centerline {$r_x(O_T(x)) = \varepsilon_0 \ \hbox{for all} \ x \in \widetilde{K_0}$.}
\vskip 0.3 cm
\noindent Afterwards, for all $x_1, x_2 \in C_{\epsilon_0}$, we have

\vskip 0.3 cm
\centerline {$\|Tx_1 - Tx_2\| \leq \frac{1}{2} (r_{x_1}(O_T(x_1)) + r_{x_2}(O_T(x_2)))$.}

\vskip 0.3 cm
\centerline {$ \leq \frac{1}{2} (\varepsilon_0 + \varepsilon_0) = \varepsilon_0$.}
\vskip 0.3 cm
\noindent This proves that

\vskip 0.3 cm
\centerline {$ \delta(T(C_{\epsilon_0})) \leq \varepsilon_0$.}
\vskip 0.3 cm
\noindent Thus

\vskip 0.3 cm
\centerline {$ \delta (\widetilde{K_0}) = \delta(\overline{co}(T(C_{\epsilon_0}))) = \delta(T(C_{\epsilon_0})) \leq \varepsilon_0$.}
\vskip 0.3 cm
\noindent Consequently,

\vskip 0.3 cm
\centerline {$ \delta (\widetilde{K_0}) = \varepsilon_0$.}
\vskip 0.3 cm
\noindent Finally,

\vskip 0.3 cm
\centerline {$ \delta (\widetilde{K_0}) = r_z(O_T(z))$ for all $z \in \widetilde{K_0}.$}
\vskip 0.3 cm
\noindent 2) Assume that $T$ has no fixed points. So, by 1), we obtain the existence of a $T$-invariant closed convex subset $D_0$ of $C$ with $\delta(D_0) > 0$ and $r_x(O_T(x)) = \delta(D_0)$ for all $x \in D_0$, which is a contradiction.

\vskip 0.3 cm
\begin{theorem}\label{th 3.6}\rm Let $X$ be a Banach space. Then the following assertions are equivalent:
\vskip 0.3 cm
\noindent $\imath)$ $X$ has weak normal structure;
\vskip 0.3 cm
\noindent $\imath \imath)$ \ Every orbitally Kannan mapping which diminishes the radius of orbits on a weakly compact convex subset of $X$ has a (unique) fixed point.
\end{theorem}
\proof $\imath) \Longrightarrow \imath \imath)$ is an immediate consequence of the assertion $1)$ in Theorem 4.5.
\vskip 0.3 cm
\noindent $\imath \imath)\Longrightarrow \imath)$ Assume the contrary that $X$ does not have weak normal structure then there exists a nontrivial weakly compact convex subset $K_0$ of $X$ such that $K_0$ does not have normal structure. Following the same notations in the proof of Theorem 3.2, we can observe that the mapping $T$ is an orbitally Kannan mapping which dimishes the radius of orbits since we have
\vskip 0.3 cm
\centerline {$ r_x(O_T(x)) = \delta(\widetilde{K_0})$ for all $x \in \widetilde{K_0}.$}
\vskip 0.3 cm
\noindent However, as indicated in Theorem 2.2 of \cite{Ami}, $T$ fails to have fixed points.

\vskip 0.3 cm
\begin{corollary}\label{cor 3.1}\rm Let $X$ be a Banach space having weak quasi-normal structure but fails to have weak normal structure. Then there exists a nontrivial weakly compact convex subset $K$ of $X$ and a mapping $T: K\longrightarrow K$ such that
\vskip 0.3 cm
\noindent $\imath)$ $T$ is not a Kannan mapping.
\vskip 0.3 cm
\noindent $\imath \imath)$ $T$ is an orbitally Kannan mapping which diminishes the radius of orbits
\vskip 0.3 cm
\noindent $\imath \imath \imath)$ $T$ has no fixed points.
\end{corollary}
\vskip 0.3 cm
\noindent As an immediate consequence, we have

\vskip 0.3 cm
\begin{corollary}\label{cor 3.2}\rm Let $X$ one of the following Banach spaces
\vskip 0.3 cm
\noindent $\imath)$ $L^{1}([0, 1]);$
\vskip 0.3 cm
\noindent $\imath \imath)$ $C([0, 1]);$
\vskip 0.3 cm
\noindent $\imath \imath \imath)$ the space $c_0$;
\vskip 0.3 cm
\noindent $\imath v)$ The space $l^{2}_{|.|}$ of Example 3.1.
\vskip 0.3 cm
\noindent Then there exists a nontrivial weakly compact convex subset $K$ of $X$ and a mapping $T: K \longrightarrow K$ such that $T$ satisfies $\imath), \imath \imath), \imath \imath \imath)$ of Corollary 4.1.
\end{corollary}

\vskip 0.3 cm
\begin{remark}\label{rem 3.1}\rm Alspach's example
\vskip 0.3 cm
\noindent Let
\vskip 0.3 cm
\centerline {$K_0 = \{f \in L^{1}([0, 1]): 0 \leq f \leq 2, \|f\|_{L^{1}([0, 1])} = 1$\}.}
\vskip 0.3 cm
\noindent Let us define $T: K_0 \longrightarrow K_0$ by
\vskip 0.3 cm

\vskip 0.3 cm
\centerline{$Tf(t)=\left\{
               \begin{array}{ll}
                 2f(2t) \wedge 2, \ \ \ & \hbox{if}\ 0\leq t \leq\displaystyle \frac{1}{2}, \\
                 (2f(2t-1)-2) \vee 0, \ \ \ & \hbox{if}\ \displaystyle \frac{1}{2}< t \leq 1.
               \end{array}
             \right.$}
\vskip 0.3 cm
\noindent In \cite{Alp}, D. Alspach proved that $T$ is a fixed point free nonexpansive mapping. More precisely, $T$ is a (nonlinear) isometry, then $T$ diminishes the radius of invariant convex subsets of $K_0$. On the other hand, since $L^{1}([0, 1])$ has the weak quasi-normal structure, by using Theorem 4.1, we deduce that $T$ is not a Kannan mapping. If $h_0 \in K_0$ is a sample starting function, for example $h_0 = \displaystyle \frac{1}{2}\chi_{[0, 1]}$, we can determine the orbit $(O_T(h_0))$. Indeed, for every integer $n \geq 1$, we have
\vskip 0.3 cm
\centerline {$T^{n}h_0 = \displaystyle \frac{1}{2}(r_n + 1)$}
\vskip 0.3 cm
\noindent where $r_n$ is the n'th Rademacher function (see \cite{Sim}). But the convex subset $co\{T^{n}h_0: n \geq 1\}$ is not an invariant subset of $K_0$. It seems not trivial to determine all orbits $O_T(h)$ for $h \in K_0$ and it will be interesting to check whether $T$ is an orbitally Kannan mapping.
\end{remark}

\vskip 0.3 cm
\begin{remark}\label{rem 3.2}\rm By using the same way of the proofs in this section together with the ideas in \cite{Rou}, we can establish the dual version of our results.
\end{remark}
\section {Final Remarks and Conclusion}
\noindent Firstly, we can easily observe that our strategy in this paper can be applied to study the existence of fixed points for $n$-orbitally Kannan mappings which diminishes the radius of orbits given as follows:

\vskip 0.3 cm
\begin{definition}\label{def 3.2}\rm (see \cite{Goh1}) Let $C$ be a nonempty subset of a Banach space $X$ and let $T: C \longrightarrow C$ be a self-mapping, $T$ is said to be $n$-orbitally Kannan mapping if
\vskip 0.3 cm
\centerline {$\|Tx - Ty\|^{n}\leq \displaystyle \frac{1}{2}[(r_x(O_T(x)))^{n} + (r_x(O_T(y)))^{n}]$ for all $x, y \in C$}
\vskip 0.3 cm
\noindent where $n$ is an integer $\geq 1$.
\end{definition}
\vskip 0.3 cm
\noindent By dropping the convexity assumption, W. A. Kirk and N. Shahzad (see \cite{Kir}) have obtained interesting results associated with orbitally wrt nonexpansive mappings. The proof of their main result can be adapted to establish a same result for the case of mappings which diminish the radius of invariant convex subsets. To be clear, we have

\vskip 0.3 cm
\begin{theorem}\label{the 3.2}\rm (see Theorem 5.2 in \cite{Kir}) Let $K$ be a weakly compact subset of a Banach space $X$. Assume that $T: K \longrightarrow K$ is a mapping which diminishes the radius of invariant weakly closed subsets of $K$ and assume that for all $x \in K$ with $x \neq Tx$,
\vskip 0.3 cm
\centerline {$\displaystyle \inf_{m \in \mathbb{N}}\{ \displaystyle \limsup_n \|T^{n}x - T^{m}x\|\} < \delta(O_T(x)).$}
\vskip 0.3 cm
\noindent Then $T$ has a fixed point.
\end{theorem}
\vskip 0.3 cm
\noindent As a consequence of the above theorem, we can derive the following

\vskip 0.3 cm
\begin{corollary}\label{cor 3.2}\rm Let $K$ be a weakly compact subset of a Banach space $X$. Assume that $T: K \longrightarrow K$ is a mapping which diminishes the radius of invariant weakly closed subsets of $K$. Assume that for every $x \in K$, there exists $\alpha(x) \in (0, 1)$ and $k \in \mathbb{N}$ such that for all integer $n \geq k$, we have
\vskip 0.3 cm
\centerline {$\|Tx - T^{n}x\| \leq \alpha(x) \|x - T^{n}x \|.$}
\vskip 0.3 cm
\noindent Then $T$ has a fixed point.
\end{corollary}
\vskip 0.3 cm
\centerline {\bf  Conclusion}
\vskip 0.3 cm
\noindent Our results refine well known contributions in the literature, we quote for example those given in (\cite {Ami, Kir, Won}). Techniques of our proofs are based essentially on some ideas observed in the indicated papers and Chapter 4 of \cite{Kha}.
\section*{Funding}

\noindent This work supported by the research team RPC (Controllability and Perturbation Results) in the laboratory of Informatics and Mathematics (LIM) at the university of Souk-Ahras (Algeria).

\section*{Competing interests}
\noindent The authors declare that they have no competing interests.

\section*{Authors contributions}
\noindent All authors contributed equally and significantly in writing this
article. All authors read and approved the manuscript.

\end{document}